\documentclass[12pt]{article}
\usepackage[english]{babel}
\usepackage[letterpaper,top=2cm,bottom=2cm,left=3cm,right=2cm,marginparwidth=1.75cm]{geometry}
\usepackage{amsfonts}
\usepackage{hyperref,amsthm,amsmath,amssymb,commath, color,multirow,graphicx}
\usepackage{inputenc}
\usepackage{wrapfig, float,array}
\usepackage{xcolor}

\begin{document}

\begin{center}
{\sc\textbf{{\Large {\ FRACTIONAL STURM–LIOUVILLE PROBLEM ON METRIC GRAPHS}}}}\\
\medskip

\, A. A. Turemuratova$^{2,b}$,\, R. Ch. Kulaev$^{3,4,c}$, Z. A. Sobirov$^{1,2,a}$
\\
\medskip
{\normalsize{\sl $^{1}$ V.I. Romanovskiy Institute of Mathematics, Uzbekistan Academy of Sciences, Tashkent, Uzbekistan\\ $^{2}$ National University of Uzbekistan, Tashkent, Uzbekistan\\
$^3$ Southern Mathematical Institute, Vladikavkaz Scientific Center, Russian Academy of Sciences\\
$^4$ North Ossetian State University after K.L.~Khetagurov

E-mail: $^{a}$z.sobirov@nuu.uz,
$^{b}$turemuratova\_a@nuu.uz,
$^{c}$kulaevrch@mail.ru}}\\
\end{center}

\begin{abstract}
In the present paper, we investigate the fractional analog of the Sturm-Liouville problem on a metric graph using a combination of left Riemann-Liouville and right Caputo fractional derivatives. This combination creates a symmetric and positive analog of the Sturm-Liouville operator. We demonstrated that the operator has a countable number of eigenvalues converging to infinity and analyzed the convergence of the series of the reciprocal eigenvalues, providing estimates for the eigenfunctions.
\end{abstract}

\textbf{2010 Mathematical Subject Classification:}{ 34A08, 34B45,  45C05}

\textbf{Keywords:} {metric graph, fractional derivatives, fractional integral, fractional Sturm– Liouville
problem, eigenvalue, eigenfunction.}

\section{Introduction}

\ \ \ \ Exploring differential equations on graphs is crucial for advancing mathematical and physical understanding, and a wealth of compelling research is available in this area, for example, \cite{Pokorny, A1,A2,A3}. In the paper \cite{Zvereva}, a system of strings positioned in equilibrium along a star-graph structure is investigated. Ashurov et al. study the initial-boundary value problem for fractional order parabolic equation on a metric star graph in Sobolev spaces in \cite{A4}.

The Sturm-Liouville (S-L) problem is a specific class of differential equations that includes additional constraints known as boundary values on their solutions. These equations are frequently encountered in both classical physics, such as thermal conduction, and quantum mechanics, like the Schr\"odinger equation. They describe processes in which an external value (the boundary value) is kept constant while the system in question transmits some form of energy.

In 1836-1837, the French mathematicians Charles-François Sturm and Joseph Liouville independently tackled the problem of heat conduction through a metal bar. In doing so, they developed techniques for solving a wide range of differential equations. The simplest form of these equations can be expressed as $[p(x)y']{'} +[q(x) - \lambda r(x)]y = 0$, where $y$ represents a physical quantity (or the quantum mechanical wave function), and $\lambda$ is a spectral parameter, that ensures $y$ meets the boundary conditions at the endpoints of the interval over which the variable $x$ is defined. If the functions $p$, $q$, and $r$ satisfy certain suitable conditions, this equation will yield a family of solutions known as eigenfunctions, each corresponding to the respective eigenvalue.


The study of spectral theory on graphs represents a significant and impactful area of research. In \cite{Zavgorodny, Pokornyi, Kulaev3} the properties (root completeness, estimation of the first eigenvalue, oscillatory properties) of the eigenvalues of the operator S-L on the graph are studied. Higher-order spectral problems are studied in works \cite{Kulaev1,Kulaev2}. The paper \cite{Bondarenko} addresses the inverse spectral problem for the matrix S-L operator (SLO) with general self-adjoint boundary conditions, aiming to provide a constructive solution. Burlutskaya et al. \cite{Burlutskaya} investigate a boundary value problem on a geometric star graph, focusing on a second-order differential equation with impulsive singularities in its coefficients and right-hand side. This is influenced by localized external loads such as elastic supports and concentrated forces. Another significant contribution to the field of spectral theory on graphs is the work \cite{Friedman} by Friedman and colleagues. They develop a wave equation for graphs that is based on a type of graph Laplacian and give some applications of this wave equation to eigenvalue/geometry inequalities on graphs. And also integro-differential equations are studied in \cite{T1,T2}.

The exploration of fractional Sturm–Liouville problems (FSLP) represents a crucial and dynamic field within the study of differential equations. Engaging with FSLP opens new avenues for research and application, highlighting its significance in advancing mathematical understanding. In \cite{Mutlu} the authors discuss the relationship between the self-adjointness of the SLO operator in the new Hilbert space and in the classical $L_2$ space, providing insights into the mathematical properties of these operators.  Klimek et al. have explored the FSLP, producing findings that are helpful for fractional order equations. The authors investigate the eigenvalue and eigenfunction properties of these fractional SLOs, highlighting that the Legendre Polynomials derived from the Fractional Legendre Equation are identical to those from the standard Legendre equation, although their eigenvalues differ \cite{klimek2013}. In \cite{klimek2016} the results indicate that under certain assumptions, there exists an infinite increasing sequence of eigenvalues, each associated with a unique continuous eigenfunction, which forms an orthogonal set of solutions. The paper successfully proves the existence of a purely discrete, countable spectrum for the FSLP under homogeneous mixed boundary conditions \cite{klimek2018}.

The paper is organized as follows: In Section 2, we provide the preliminaries. In Section 3, we consider the FSLP for metric graphs. Some properties of eigenvalues and eigenfunctions are discussed in Section 4. Finally, we summarize in Section 5.


\section{Preliminaries}

\ \ \ \ In this section, we present essential definitions and properties.

\textbf{Definition 1.} \cite{Kilbas} The left and right Riemann–Liouville fractional integrals of order $0<\alpha<1$  for a function $y(x)\in L_1(0,l)$ are, respectively, defined by
\begin{equation} \nonumber
I^{\alpha}_{0,x}y(x)\equiv{\mathrm{D}}^{-\alpha}_{0,x}y(x):=\frac{1}{\Gamma(\alpha)}\int_0^x\frac{y(\xi)}{(x-\xi)^{1-\alpha}}d\xi,
\end{equation}
\begin{equation} \nonumber
I^{\alpha}_{x,l}y(x)\equiv{\mathrm{D}}^{-\alpha}_{x,l}y(x):=\frac{1}{\Gamma(\alpha)}\int_x^l\frac{y(\xi)}{(\xi-x)^{1-\alpha}}d\xi.
\end{equation}

\textbf{Lemma 1.} \cite{Kilbas} If $\alpha>0$ and $1\leq p\leq\infty$, then $I^{\alpha}_{0,x}$ and $I^{\alpha}_{x,l}$ are continuous from $L_p(0,l)$ into itself and
$$\rVert I^{\alpha}_{0,x}y \rVert_{L_p(0,l)} \leq \frac{l^\alpha}{\Gamma (\alpha+1)} \rVert y \rVert_{L_p(0,l)},\,\,\,\,\,\,\rVert I^{\alpha}_{x,l}y \rVert_{L_p(0,l)} \leq \frac{l^\alpha}{\Gamma (\alpha+1)} \rVert y \rVert_{L_p(0,l)},$$
for all $y\in L_p(0,l)$.

\textbf{Definition 2.}
\cite{Kilbas} The left and right Riemann-Liouville fractional derivatives of order $0<\alpha<1$ for a function $y$ on $[0,l]$ are respectively, defined by

$$ {\mathrm{D}}^{\alpha}_{0,x}y(x)=\frac{1}{\Gamma(1-\alpha)}\frac{d}{dx}\int_0^x\frac{y(\xi)}{(x-\xi)^{\alpha}}d\xi,$$
$${\mathrm{D}}^{\alpha}_{x,l}y(x)=\frac{-1}{\Gamma(1-\alpha)}\frac{d}{dx}\int_x^l\frac{y(\xi)}{(\xi-x)^{\alpha}}d\xi,$$
provided that, the integrals, in the right-hand sides of these expressions, are exist.

\textbf{Definition 3.}
\cite{Kilbas} The left and right Caputo fractional derivatives of order $0<\alpha<1$ for a function $y$ on $[0,l]$ are respectively, defined by

$${\partial}^{\alpha}_{0,x}y(x)=\frac{1}{\Gamma(1-\alpha)}\int_0^x\frac{y^{'}(\xi)}{(x-\xi)^{\alpha}}d\xi,
\qquad\ {\partial}^{\alpha}_{x,l}y(x)=\frac{-1}{\Gamma(1-\alpha)}\int_x^l\frac{y^{'}(\xi)}{(\xi-x)^{\alpha}}d\xi,$$
provided that, the integrals, in the right-hand sides of these expressions, are exist.

\textbf{Definition 4.} Suppose $u,v\in L_1(0,l)$ and $1/2<\beta<1$. We say that $v$ is the weak left fractional derivative of order $\beta$ of $u$, written $D^{\beta}_{0,x}u=v$, provided
$$\int\limits_0^l u(x) \partial^{\beta}_{x,l}\varphi(x) dx=\int\limits_0^l v(x)\varphi(x)dx$$
for all test functions $\varphi\in C_{c}^{\infty}(0,l)$.

We established some functional spaces.
Let $0<\beta<1$. We define $AC^{\beta,2}_{0}(0,l)=AC^{\beta,2}_{0}\left([0,l],\mathbb R\right)$ as the set of all functions that have representation
$$f(x)=\frac{c_0}{\Gamma(\beta)}x^{\beta-1}+I^{\beta}_{0,x}\phi(x)\qquad for\quad a.e.\quad x\in[0,l],$$
with some $\phi \in L_2(0,l)$, $c_0\in\mathbb R$, and by $AC_l^{\beta,2}(0,l)=AC_l^{\beta,2}\left([0,l], \mathbb R\right)$ we mean the set of all functions that have the representation
$$g(x)=\frac{d_0}{\Gamma(\beta)}(l-x)^{\beta-1}+I^{\beta}_{x,l}\psi(x)\qquad for\quad a.e.\quad x\in[0,l],$$
with some $\psi \in L_2(0,l)$,  $d_0\in \mathbb R$.

We set \cite{Idczak},\cite{9}
$\mathcal{H}_+^{\beta}(0,l)=\left\{u\in L_2(0,l), \,\, D_{0,x}^{\beta}u\in L_2(0,l) \right\}$,
$$\mathcal{H}^{\beta}_+(0,l)=AC^{\beta,2}_0(0,l) \cap L_2(0,l),$$
$$\mathcal{H}^{\beta}_{-}(0,l)=AC^{\beta,2}_l(0,l) \cap L_2(0,l).$$
According to Remark 21 in \cite{Idczak} $\mathcal{H}^{\beta}_+(0,l)=AC_0^\beta(0,l)$, for $\beta>\frac{1}{2}$.
Then it follows from the definitions of $AC^{\beta,2}_0(0,l)$ and $AC^{\beta,2}_l(0,l)$ that,
$$\rho \in \mathcal{H}^{\beta}_+(0,l) \Longleftrightarrow \rho \in L_2(0,l)\,\, and\,\, {\mathrm{D}}^{\beta}_{0,x} \rho \in L_2(0,l),$$
$$\rho \in \mathcal{H}^{\beta}_{-}(0,l) \Longleftrightarrow \rho \in L_2(0,l)\,\, and\,\, {\mathrm{D}}^{\beta}_{x,l} \rho \in L_2(0,l).$$
Furthermore, the space $\mathcal{H}^{\beta}_+(0,l)$,  equipped with the scalar product
$$(\varphi,\psi)_{\mathcal{H}^{\beta}_+(0,l)}=(\varphi,\psi)_{L_2(0,l)}+\left(\gamma(x)\mathrm{D}^{\beta}_{0,x}\varphi, \mathrm{D}^{\beta}_{0,x}\psi\right)_{L_2(0,l)},$$
and the space $\mathcal{H}^{\beta}_-(0,l)$,  equipped with the scalar product
$$(\varphi,\psi)_{\mathcal{H}^{\beta}_{-}(0,l)}=(\varphi,\psi)_{L_2(0,l)}+\left(\gamma(x)\mathrm{D}^{\beta}_{x,l}\varphi, \mathrm{D}^{\beta}_{x,l}\psi\right)_{L_2(0,l)},$$
form Hilbert spaces (see \cite{Idczak}).

\textbf{Lemma 2.}
\cite{Leugering} Let $f\in \mathcal{F}^{\beta}:=\left\{f\in C[0,l]:{\mathrm{D}}^{\beta}_{0,x}f\in L_2(0,l)\right\}$ and \\
 $g\in Q^{\beta}:=\left\{g\in C[0,l]: {\partial}^{\beta}_{x,l}g \in L_2(0,l)\right\}$, then the following holds:

$$\int_0^l{\mathrm{D}}^{\beta}_{0,x}f(x)g(x)dx=\int_0^l f(x) {\partial}^{\beta}_{x,l}g(x)dx+\big [g(x)I^{1-\beta}_{0,x}f(x) \big]^{x=l}_{x=0},$$
$$\int_0^l{\partial}^{\beta}_{x,l}g(x)f(x)dx=\int_0^l g(x) {\mathrm{D}}^{\beta}_{0,x}f(x)dx-\big [g(x)I^{1-\beta}_{0,x}f(x) \big]^{x=l}_{x=0}.$$

\textbf{Lemma 3}. Let $f\in {{L}_{2}}(0,l),$ $0<\alpha <\frac{1}{2}$. There exists the sequence of functions $\left\{ {{f}_{k}} \right\}_{k=1}^{\infty }$ such that ${{f}_{k}}\in {{C}^{\infty }}[0,l]$, $\left( I_{0,x}^{\alpha }{{f}_{k}} \right)(0)=\left( I_{0,x}^{\alpha }{{f}_{k}} \right)(l)=0,\ k=1,2,...,$ and $||{{f}_{k}}-f|{{|}_{{{L}_{2}}(0,l)}}\to 0$ as $k\to +\infty $.

\textbf{Proof.}
It is known, that ${{C}^{\infty }}[0,l]$ is dense in ${{L}_{2}}(0,l)$. So, there exist $\left\{ {{v}_{k}} \right\}\subset {{C}^{\infty }}[0,l]$, such that ${{v}_{k}}\to f$ in ${{L}_{2}}(0,l)$.
According to the following inequality
	\[\left| I_{0,t}^{\alpha }{{v}_{k}}(x) \right|=\frac{1}{\Gamma (\alpha )}\left| \int\limits_{0}^{x}{\frac{{{v}_{k}}(\xi )d\xi }{{{(x-\xi )}^{1-\alpha }}}} \right|\le \frac{1}{\Gamma (\alpha )}\underset{\xi \in [0,l]}{\mathop{\max }}\,\left| {{v}_{k}}(\xi ) \right|\int\limits_{0}^{x}{\frac{d\xi }{{{(x-\xi )}^{1-\alpha }}}}=\]
	\[=\frac{1}{\Gamma (\alpha +1)}\underset{\xi \in [0,l]}{\mathop{\max }}\,\left| {{v}_{k}}(\xi ) \right|\,\,{{x}^{\alpha }},\]
we have $\underset{x\to 0}{\mathop{\lim }}\,I_{0,t}^{\alpha }{{v}_{k}}(x)=0.$

We put \[{{\left. I_{0,x}^{\alpha }{{v}_{k}} \right|}_{x=l}}={{c}_{k}}.\] If ${{c}_{k}}=0$ for all $k\in N$, the sequence is satisfying the conditions of the theorem.

Let $\{{c}_{k}\}\ne 0$. In this case, we need to construct auxiliary functions to adjust the values of the functions $v_k$ at the right end of the interval.  For this purpose we put ${{z}_{n}}={{a}_{n}}{{x}^{n}}$, where ${{a}_{n}}=\frac{\Gamma (n+1+\alpha )}{{{l}^{n+\alpha }}\,\Gamma (n+1)}.$ Using the direct computations, one have \[{{\left. I_{0,x}^{\alpha }{{z}_{n}} \right|}_{x=l}}=1\] and
	\[||{{z}_{n}}|{{|}_{{{L}_{2}}(0.l)}}=\frac{\Gamma (n+1+\alpha )}{\Gamma (n+1)\sqrt{2n+1}}{{l}^{1-2\alpha }}\sim const\cdot \frac{1}{n{{\,}^{\frac{1}{2}-\alpha }}},\ n\to +\infty .\]
Next we choose subsequence ${{n}_{k}}$, such that $\frac{c_{k}^{\frac{2}{1-2\alpha }}}{{{n}_{k}}}\to 0$ as  $k\to +\infty .$ Then the sequence $\left\{ {{c}_{k}}{{z}_{{{n}_{k}}}} \right\}_{k=1}^{\infty }$ tends to zero in ${{L}_{2}}(0,l)$.

Now we construct the sequence $\left\{ {{f}_{k}} \right\}:\ \ {{f}_{k}}={{v}_{k}}-{{c}_{k}}{{z}_{{{n}_{k}}}}$. It is easy to see, that this sequence satisfies the conditions of the theorem. This completes the proof.

Let us define a graph and a function defined on a graph.

A graph $\mathcal{G}$ consists of a finite set of vertices $V=\left\{ {\nu_{i}} \right\}_{i=1}^{j}$ and a set $E=\left\{ {e_i} \right\}_{i=1}^{n}$ of edges. A metric graph is a graph in which each edge is associated with a positive number, called the length of this edge. Each edge $e_i$ is assigned the interval $\left( 0,{l_i} \right),\,i=\overline{1,n}$, and the coordinates $x_i$ are defined on each edge. We will say that a vertex $\nu$ is in incident with an edge $e_i$ if it is the end of this edge, and denote this as $e_i\sim \nu$. The number of elements of the set $\{e: e \sim \nu, e \in E\}$ is called vertex valency $\nu$. If the valency of a vertex is equal to one, then it is called a boundary vertex. Let $\partial \mathcal{G} \subset V$ be the set of boundary vertices of the graph.

During this work, we suppose that $\mathcal{G}$ is a connected graph with a non-empty set of boundary vertices.

For the function $u:\mathcal{G}\to R$, defined on the graph, we put $u\rvert_{e_i}=u^{(i)}$. For the functions defined on the graph, we also use vector-type notations $u=(u^{(1)},...,u^{(n)})$, $u_x=\left(\frac{\partial u^{(1)}}{\partial x},..., \frac{\partial u^{(n)}}{\partial x}\right)$,  $\int\limits_{\mathcal{G}}u d\mathcal{G}= \sum\limits_{i=1}^n \int\limits_{0}^{l_i} u^{(i)} dx$.
For $u:\mathcal{G}\to R$, $v:\mathcal{G}\to R$ we put $uv=(u^{(1)}v^{(1)},u^{(2)}v^{(2)},\ldots,u^{(n)}v^{(n)})$.

We set  $C[\mathcal{G}]$ as a space of functions $u$ defined on the graph and uniformly continuous on each of its edges.
The space $C(\mathcal{G})$ of functions $u:\mathcal{G}\rightarrow \mathbb{R}$ defined on the entire graph $\mathcal{G}$ and continuous on the graph.

\textbf{Path integral.}
\begin{figure}[ht!]
\begin{center}
\includegraphics[width=60mm]{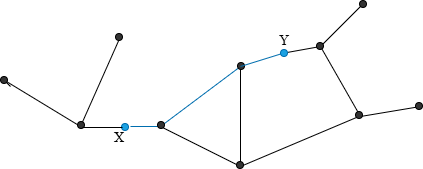}
\end{center}
\caption{Metric graph}
\label{fig1}
\end{figure}
Let $A,B\in \mathcal{G}$ (Fig. 1). By $\gamma_{AB}$ we denote the shortest path on the graph, which connects the points $A$ and $B$. We put $\rho(A,B)=\rvert \gamma_{AB}\rvert$, where $\rvert \gamma_{AB}\rvert$ is a length of the path $\gamma_{AB}$.

Consider an isometric map
$$F:\gamma_{AB}\rightarrow\left[0,\rvert \gamma_{AB}\rvert\right], \quad F(X)=\rho(A,X).$$
Then we define path integral by
$$\int_{\gamma_{AB}}f(x)dx=\int\limits_0^{\rvert \gamma_{AB}\rvert}\tilde{f}(\xi)d\xi,$$
where $f(F(x))=f(x)$. It is clear, that if $f\in H^1(\mathcal{G})=\left(\bigoplus_{i=1}^{n} H^1(e_i)\bigcap C(\mathcal{G})\right)$, then
$$\int_{\gamma_{AB}}f_xdx=f(B)-f(A), \eqno(1)$$
for any $A,B\in \mathcal{G}$.

By $\int fd\gamma_{AB}$ we denote the integral over the subgraph $\gamma_{AB}$ defined above. One should notice that these two integrals are different:
$$\int_{\gamma_{AB}}f(x)dx=-\int_{\gamma_{BA}}f(x)dx,$$
while $\int f(x)d\gamma_{AB}=\int f(x)d\gamma_{BA}$.

Let $\mathcal{G}$ is a connected graph, $\partial \mathcal{G}\ne \emptyset$, $f:\mathcal{G}\rightarrow \mathbb{R}$, $f\in H^1(\mathcal{G})$ and $f(\nu)=0$ for some $\nu \in \partial\mathcal{G}$. Then, from (1) it follows that
$$f(X)=\int_{\gamma_{\nu X}} f_{x}(x)dx. \eqno(2)$$
From (2) it is easy to get
$$\int_{\mathcal{G}}f^2(X)d\mathcal{G}_X\leq\int_{\mathcal{G}}\left(\int_{\gamma_{\nu X}}f_xdx\right)^2d\mathcal{G}_X\leq \int_{\mathcal{G}}\left(\int_{\mathcal{G}}\rvert f_x(x)\rvert d\mathcal{G}\right)^2d\mathcal{G}\leq \rvert \mathcal{G}\rvert^2\int_{\mathcal{G}}f_{x}^2d\mathcal{G}, \eqno(3)$$
where $\rvert \mathcal{G} \rvert$ is a total length of the graph.

\textbf{Definition 5.}  The space ${L_2}(\mathcal{G} )$ on $\mathcal{G} $ consists of functions that are measurable and square-integrable on each edge ${e_i},\,\,i=\overline{1,n}$ with the scalar product and the norm:
$${\left( u\left( x \right),v\left( x \right) \right)}_{L_2(\mathcal{G} )}=\int\limits_{\mathcal{G} }{u\cdot v\,d\mathcal{G} }=\sum\limits_{i=1}^{n}{\int\limits_{0}^{l_i}u^{(i)}(x)\cdot v^{(i)}(x) dx},  $$

$$\left\| u \right\|_{L_2\left( \mathcal{G}  \right)}^{2}=\sum\limits_{i}{\left\| u \right\|^2_{L_2(e_i)}}.$$

In other words, $L_2(\mathcal{G})$ is the orthogonal direct sum of spaces $L_2(e_i),\,\, i=\overline{1,n}$.\\

We put the space $\mathcal{H_+}^{\beta}(\mathcal{G})=\{u \in \bigoplus_{i=1}^{n} \mathcal{H}^{\beta}_{+}(0,l_i),\,\, I^{1-\beta}_{0,x}u\in C(\mathcal{G})\}$  equipped
with the norm given by
$$\rVert u \rVert^2_{\mathcal{H}_+^{\beta}(\mathcal{G})}=\sum_{i=1}^n \|u^{(i)}\|^2_{\mathcal{H}^{\beta}_{+}(0,l_i)}.$$

\section{Main problem}

\ \ \ \ The purpose of this paper is to present existence results for space-time fractional diffusion equations by employing the method of separation of variables. This approach relies significantly on fractional Sturm-Liouville theory, specifically on the existence of eigenvalues and their corresponding eigenfunctions for the following fractional differential equation on the graph $\mathcal{G}$
$$
\mathcal{L}u=\lambda r u,\quad x\in e_i, \quad \frac{1}{2}<\beta<1, \eqno(4)
$$
where $u:\mathcal{G}\rightarrow\mathbb{R}$, $r:\mathcal{G}\rightarrow\mathbb{R}$, and
$$\mathcal{L}u\big \rvert_{e_i}=\partial_{x,l_i}^\beta\left(p^{(i)}(x)D_{0,x}^\beta u^{(i)}(x)\right), \eqno(5)$$
where $0<a_0\leq p^{(i)}\leq a_1<+\infty$, $0<R_0\leq r^{(i)}\leq R_1<+\infty$ and we suppose that $r, p\in C[\mathcal{G}]$.
At each interior vertex $\nu$, continuity conditions and transmission conditions
$$
\begin{cases}
I_{0,x}^{1-\beta}u^{(i)}(x) \,\,are\,\, continuous\,\, in\,\, \nu, \\
 \sum\limits_{e_i\sim \nu}\sigma_{e_i,\nu} {p^{(i)}(\nu)D_{0,x}^{\beta}u^{(i)}(\nu)}=0, \quad \nu\in V\setminus  \partial \mathcal{G},
\end{cases} \eqno(6)
$$
are specified, where $\sigma_{e_i,\nu}=1$  if $\nu$ is the right end of the edge $e_i$, $\sigma_{e_i,\nu}=-1$ if $\nu$ is the left end of the edge $e_i$, and at each boundary vertex, Dirichlet conditions
$$I_{0,x}^{1-\beta}u(\nu)=0,\quad \nu\in\partial \mathcal{G}, \eqno(7)$$
are specified.

We define the following operators on the graph
$$D_{+}^{\beta}u=\left\{D_{0,x}^{\beta}u^{(1)},D_{0,x}^{\beta}u^{(2)},\ldots,D_{0,x}^{\beta}u^{(n)}\right\},\,\,\,\,D_{-}^{\beta}u=\left\{D_{x,l_1}^{\beta}u^{(1)},D_{x,l_2}^{\beta}u^{(2)},\ldots,D_{x,l_n}^{\beta}u^{(n)}\right\},$$
which can be considered as a fractional order of the gradient operator. The operators $\partial_{\pm}^\beta$ can be defined similarly.

We put $\mathcal{H}^\beta_{+,Dir}(\mathcal{G})=\left\{u\in \mathcal{H}_{+}^\beta(\mathcal{G}): \left.I^{1-\beta}_{0,x}u\right|_{\partial\mathcal{G}}=0 \right\}$.

\textbf{Definition 6.} The function $u(x)\in \mathcal{H}^{\beta}_{+,Dir}(\mathcal{G})$ is a generalized solution to problem (4)-(7) if the equality
$$\left(p(x)D_{+}^\beta u(x), D_{+}^\beta w(x)\right)_{L_2(\mathcal{G})}=\lambda \left(r(x)u(x),w(x)\right)_{L_2(\mathcal{G})}\eqno(8)$$
holds for any $w(x)\in \mathcal{H}^{\beta}_{+,Dir}(\mathcal{G})$.

\textbf{Remark 1}. For smooth enough functions, the identity (8) can be obtained from equation (4) by multiplying $w$ and integrating by part.

For further investigation, we need the following equivalence of norms.

\textbf{Lemma 4.} In $\mathcal{H}^{\beta}_{+,Dir}(\mathcal{G})$ the norms $\| u \|_{1,\mathcal{H}_+^{\beta}(\mathcal{G})}=\left(\int\limits_{\mathcal{G}} \left(\sqrt{p(x)} D^\beta_+u\right)^{2}d\mathcal{G}\right)^{1/2}$  and $\|\cdot\|_{\mathcal{H}_+^\beta(\mathcal{G})}$  are equivalent.
\begin{proof}
According to Theorem 24 in \cite{Idczak} we have the following equivalence of norms
	\[{{C}_{1}}\left( |I_{0,x}^{1-\beta }{{u}^{(j)}}(0){{|}^{2}}+||D_{0,x}^{\beta }{{u}^{(j)}}||_{{{L}_{2}}(0,{{l}_{j}})}^{2} \right)\le ||{{u}^{(j)}}||_{\mathcal{H}_{+}^{\beta }(0,{{l}_{j}})}^{2}\le {{C}_{2}}\left( |I_{0,x}^{1-\beta }{{u}^{(j)}}(0){{|}^{2}}+||D_{0,x}^{\beta }{{u}^{(j)}}||_{{{L}_{2}}(0,{{l}_{j}})}^{2} \right).\]

Using this inequality, we get
	\[||u||_{\mathcal{H}_{+}^{\beta }(\mathcal{G})}^{2}\le {{C}_{2}}\left( \sum\limits_{j=1}^{n}{|I_{0,x}^{1-\beta }{{u}^{(j)}}(0){{|}^{2}}}+||D_{0,x}^{\beta }u||_{{{L}_{2}}(\mathcal{G})}^{2} \right).\]	

Now, using the inequality (2) for $f=I_{0,x}^{1-\beta }u$, we have
	\[|I_{0,x}^{1-\beta }u(X)|=\left| \int\limits_{{{\gamma }_{vX}}}{D_{0,x}^{\beta }u(x)dx} \right|\le \int\limits_{\mathcal{G}}{\left| D_{0,x}^{\beta }u(x) \right|d\mathcal{G}}\le \sqrt{|\mathcal{G}|}||D_{0,x}^{\beta }u(x)|{{|}_{{{L}_{2}}(\mathcal{G})}}.\]	

Combining the last two inequalities we get
	\[||u||_{\mathcal{H}_{+}^{\beta }(\mathcal{G} )}^{2}\le {{C}_{2}}\left( n|G|+1 \right)||D_{0,x}^{\beta }u||_{{{L}_{2}}(\mathcal{G})}^{2}.\]

On the other hand,
	\[||u||_{\mathcal{H}_{+}^{\beta }(\mathcal{G} )}^{2}=||u||_{{{L}_{2}}(G)}^{2}+||D_{0,x}^{\beta }u||_{{{L}_{2}}(\mathcal{G} )}^{2}\le ||D_{0,x}^{\beta }u||_{{{L}_{2}}(\mathcal{G})}^{2}.\]

As $0<a_0\leq p^{(i)}\leq a_1<+\infty$ we have $\sqrt{a_0}||u||_{\mathcal{H}_{+}^{\beta }(\mathcal{G} )}\leq \rVert u \rVert_{1,\mathcal{H}_+^{\beta}(\mathcal{G})}\leq \sqrt{a_1}||u||_{\mathcal{H}_{+}^{\beta }(\mathcal{G})}$.

This completes the proof of the lemma.
\end{proof}

Due to $0<R_0\leq r^{(i)}\leq R_1<+\infty$, we could also introduce an equivalent norm in $L_2(\mathcal{G})$
$$\|u\|_{1,L_2(\mathcal{G})}^2=\int_{\mathcal{G}}ru^2d\mathcal{G}.$$

\textbf{Theorem 1.} For the eigenfunctions $u_k$ and the eigenvalues $\lambda_k$ of problem (4)-(7), the following statements are true:

1) $ 0<\lambda_1\leq\lambda_2\leq\lambda_3\leq \ldots $,

2) $\lambda_k \rightarrow \infty$,

3) the $u_k$ form a complete orthonormal basis for $L_2(\mathcal{G})$ with respect to norm $\|\cdot\|_{1,L_2(\mathcal{G})}$.

\textbf{Proof.}
Consider the Rayleigh quotient
$$\mathcal{R}(u)=\frac{\int_{\mathcal{G}} p(x)(D_{+}^\beta u(x))^2 d\mathcal{G}}{\int_{\mathcal{G}} r u^2 d\mathcal{G}}$$
which is certainly defined for $u\in \mathcal{H}_+^\beta (\mathcal{G})$. Let us put
$$\lambda=\underset{\underset{\rVert u \rVert \ne 0}{\mathop{u\in \mathcal{H}^{\beta}_{+,Dir}(\mathcal{G})}}}{\mathop{inf}}\mathcal{R}(u)=\underset{\underset{\rVert u \rVert =1}{\mathop{u\in \mathcal{H}^{\beta}_{+,Dir}(\mathcal{G})}}}{\mathop{inf}}\mathcal{R}(u). \eqno(9)$$
Let $f_1,f_2,\ldots$ be a minimizing sequence for $\mathcal{R}$ in $\mathcal{H}^{\beta}_{+,Dir}(\mathcal{G})$ with
$$\mathcal{R}(f_j)\rightarrow \lambda.$$
We may assume $\int_{\mathcal{G}}r{f_j}^2d\mathcal{G}=1$. Then $\rVert f_j \rVert_{\mathcal{H}_+^\beta}\leq \sqrt{1+\lambda}.$ Hence by passing to a subsequence we may assume that the $f_j$ converge weakly in $\mathcal{H}_+^\beta$ to a $f\in \mathcal{H}^\beta_{+,Dir}(\mathcal{G}).$
According to Corollary 32 in \cite{Idczak}, the embedding $\mathcal{H}_+^\beta(\mathcal{G})\subset L_2(\mathcal{G})$ is compact. So, $\int_{\mathcal{G}}rf^2d\mathcal{G}=1$.

By the weak convergence in $\mathcal{H}_+^\beta (\mathcal{G})$ we have
$$\mathcal{R}(f)\leq \lim inf\mathcal{R}(f_j)=\lambda.$$
On the other hand $\mathcal{R}(f)\geq \lambda$. It follows from here $\mathcal{R}(f)=\lambda$. This means that $f$ minimizes $\mathcal{R}$ over all of $\mathcal{H}^{\beta}_{+,Dir}(\mathcal{G})$.

Now we claim that $f$ is our desired eigenfunction and $\lambda=\mathcal{R}(f)$ its eigenvalue. This is seen by setting
$$\mathcal{R}(f+tg)\geq \lambda,\ \forall g \in \mathcal{H}^{\beta}_{+,Dir}(\mathcal{G}),\ \forall t \in \mathbb{R}.$$
$$\mathcal{R}(f+tg)=\frac{\rVert f+tg \rVert^2_{\mathcal{H}^{\beta}_{+,Dir}(\mathcal{G})}}{\rVert f+tg \rVert^2_{L_2(\mathcal{G})}}\geq \lambda$$
$$\rVert f+tg \rVert^2_{\mathcal{H}^{\beta}_{+,Dir}(\mathcal{G})} \geq \lambda \rVert f+tg \rVert^2_{L_2(\mathcal{G})} $$
$$\rVert f\rVert^2_{\mathcal{H}^{\beta}_{+,Dir}(\mathcal{G})}+2t\left( f,g\right)_{\mathcal{H}^{\beta}_{+,Dir}(\mathcal{G})}+t^2\rVert g\rVert^2_{\mathcal{H}^{\beta}_{+,Dir}(\mathcal{G})}\geq \lambda \rVert f \rVert^2_{L_2(\mathcal{G})}+2\lambda t\left( f,g\right)_{L_2(\mathcal{G})}+\lambda t^2\rVert g\rVert^2_{L_2(\mathcal{G})}$$

Taking into account that
$$ \rVert f\rVert^2_{\mathcal{H}^{\beta}_{+,Dir}(\mathcal{G})}-\lambda\rVert f \rVert^2_{L_2(\mathcal{G})}=0,$$
we have
$$At^2+2Bt\geq 0,\ \mbox{ for all  } t\in \mathbb{R},$$
where $A=\rVert g\rVert^2_{\mathcal{H}^{\beta}_{+,Dir}(\mathcal{G})}-\lambda\rVert g \rVert^2_{L_2(\mathcal{G})}\geq 0$, $B=\left( f,g\right)_{\mathcal{H}^{\beta}_{+,Dir}(\mathcal{G})}-\lambda\left( f,g\right)_{L_2(\mathcal{G})}$. From this $B=0$ or  $\left( f,g\right)_{\mathcal{H}^{\beta}_{+,Dir}(\mathcal{G})}-\lambda\left( f,g\right)_{L_2(\mathcal{G})}=0$, $\forall g\in \mathcal{H}^{\beta}_{+,Dir}(\mathcal{G})$. This proves our claim.

Start by setting $u_1 = f$ and $\lambda_1 =\lambda$. Then, repeat the same process, but this time minimize $\mathcal{R}$ over functions that are orthogonal to $u_1$. Using the same process, we find another eigenpair $(u_2,\lambda_2)$ with $\lambda_1\leq \lambda_2$ and $u_1$ orthogonal to $u_2$. Repeat this process again, minimizing over functions that are orthogonal to both $u_1$ and $u_2$. By doing this, we obtain a sequence of orthogonal eigenpairs $(u_k,\lambda_k)$ with $\lambda_k$ increasing as $k$ increases.

Next, we show that $\lambda_k\rightarrow \infty$. We consider the opposite assumption, which is that $\lambda_k\leq a$. Previously, we assumed that the following conditions are satisfied for eigenfunctions $u_k$
$$\rVert u_k\rVert_{L_2(\Gamma)}=1,\quad \rVert u_k\rVert_{\mathcal{H}^{\beta}(\Gamma)}\leq \sqrt{1+a}.$$
So, the sequence $u_k$ is bounded. Therefore, we can obtain a subsequence $u_{k_j}$ that weakly converges to a limit $u$ in $\mathcal{H}_+^{\beta}(\Gamma)$. But since the $u_{k_j}$ are
orthogonal, $u$ would be orthogonal to all $u_{k_j}$ and therefore to itself,
$$0=(u_{k},u_{k_j})=(u_{k},u).$$
As we take the limit as $k$ approaches infinity, we discover that $(u, u) = 0$ leads us directly to the conclusion that $u = 0$. This strongly indicates that $u_{k_j}$ converges weakly to $0$. This challenges our initial assumption that the values $\lambda_k$ remain bounded. Therefore, we are led to conclude that $\lambda_k\rightarrow \infty$.

 To demonstrate that the functions $u_k$ are complete, we start by assuming that they are not. If the $u_k$ are not complete, then there exists a nonzero function $v \in L_2(\mathcal{G})$ that is orthogonal to all $u_k$. According to Lemma 3, for any $\varepsilon > 0$  we can construct a function $v_{\varepsilon} \in \mathcal{H}^{\beta}_{+,Dir}[\mathcal{G}]\subset \mathcal{H}^{\beta}_{+,Dir}(\mathcal{G})$ such that $ \| v - v_{\varepsilon} \|_{1,L_2(\mathcal{G})} < \varepsilon $. For sufficiently small $\varepsilon$, the function $h$, which represents the projection of $v_{\varepsilon}$ onto the complement of the span of $u_k$'s, will have the following properties: (1) it is non-zero, (2) it is orthogonal to all $u_k$, and (3) it belongs to $\mathcal{H}^{\beta}_{+,Dir}(\mathcal{G})$.
Given that the Rayleigh quotient $ \mathcal{R}(h)$ must provide an upper bound for the $\lambda_k$ values due to their minimizing property, we conclude that the $\lambda_k$ are bounded. However, this conclusion leads to a contradiction since we already know that such a boundedness is impossible.

The proof of the theorem is complete.

\textbf{Remark 2.}
 We notice, that the system of functions $\{u_k\}$ is complete in $\mathcal{H}^{\beta}_{+,Dir}(\mathcal{G})$.

Indeed, if the function $w\in \mathcal{H}^{\beta}_{+,Dir}(\mathcal{G})$ is orthogonal to all functions $u_k$ in $\mathcal{H}^{\beta}_{+,Dir}(\mathcal{G})$, i.e.,  $\left(p(x)D_{+}^\beta u(x), D_{+}^\beta w(x)\right)_{L_2(\mathcal{G})}=0$, then from (8) it follows $ \left(r(x)u_k(x),w(x)\right)_{L_2(\mathcal{G})}=0$ for  $k=1,2,3,...$. So, $w=0$.

\section{Some properties of eigenvalues and eigenfunctions}
\ \ \ \ Here, we first consider the following  fractional equation in the interval $(0,l)$
$$\partial_{x,l}^\beta\left(p(x)D_{0,x}^\beta u(x)\right)=\lambda r(x)u(x), \eqno(10)$$
with boundary conditions
$$I^{1-\beta}_{0,x}u(0)=I^{1-\beta}_{0,x}u(l)=0. \eqno(11)$$
This problem is a particular case of the Sturm-Liouville problem considered in Section 3.

Let us denote by $\mathcal{T}$ the right inverse operator of $\mathcal{L}_0$. Specifically,
 $$\mathcal{L}_0u(x)=\partial_{x,l}^\beta\left(p(x)D_{0,x}^\beta u(x)\right)=f(x),$$
 $$u(x)=I_{0,x}^{\beta}\left(\frac{1}{p(x)}I_{x,l}^{\beta}f(x)\right)-\frac{\int\limits_{0}^{l}\frac{1}{p(x)}I_{x,l}^{\beta}f(x)dx}{\int\limits_{0}^{l}\frac{1}{p(x)}dx}I_{0,x}^{\beta}\left(\frac{1}{p(x)}\right)=\mathcal{T}f.$$

 Observe that operator $\mathcal{T}$ can be expressed as
an integral operator with kernel $\mathcal{K} = \mathcal{K}_1+ \mathcal{K}_2$, namely
$$(\mathcal{T}f)(x)=\int\limits_{0}^{l}\mathcal{K}(x,s)f(s)ds,$$
where the explicit form of the parts of the kernel looks as follows
$$\mathcal{K}_1(x,s)=\int\limits_0^{min\{x,s\}}\frac{(x-t)^{\beta-1}}{\Gamma(\beta)}\cdot\frac{(s-t)^{\beta-1}}{\Gamma(\beta)}\cdot\frac{1}{p(t)}dt,$$
$$\mathcal{K}_2(x,s)=-\frac{1}{\left.I_{0,x}^{\beta}\frac{1}{p(x)}\right|_{x=l}}\cdot I_{0,x}^{\beta}\frac{1}{p(x)}\cdot I_{0,s}^{\beta}\frac{1}{p(s)}.$$

In the case of the interval, the eigenvalues and eigenfunctions have the following property.

\textbf{Proposition 1.}
    $\sum\limits_{k=1}^{\infty}\frac{1}{\lambda_k}$ is convergent.
\begin{proof}
   It is easy to notice that if $\lambda_k$ is eigenvalue of the operator $\mathcal{L}_0$ with domain $H^{\beta}_{+,Dir}(0,l)$, then $\frac{1}{\lambda_k}$ be eigenvalue of the integral operator $\mathcal{T}f=\int\limits_0^lK(x,s)f(s)ds$.
According to \cite{klimek2016}, the function $\mathcal{K}_1(x,s)$ is continuous on $[0,l]\times[0,l]$. The function $\frac{1}{p(x)}$ is continuous on $[0,l]$, so, $\mathcal{K}_2(x,s)$ is also continuous on $[0,l]\times [0,l]$. So, we can use the trace formula for integral operator $\mathcal{T}$ to get
$$\sum\limits_{k=1}^{\infty}\frac{1}{\lambda_k}=\int\limits_0^lK(x,x)dx<+\infty.$$
\end{proof}

We notice that a similar problem for the operator $\tilde{L}_0=\frac{1}{\omega_\alpha(x)}\partial^{\alpha}_{x,l}p(x)\partial^{\alpha}_{0,x}$ is considered in \cite{klimek2016}. In our case, we improved the results in \cite{klimek2016}, where they obtained convergence of $\sum\limits_{k=1}^\infty\frac{1}{\lambda_k^2}$.

\textbf{Proposition 2.}
$$ \quad \rvert I_{0,x}^{1-\beta}u_k\rvert\leq\sqrt{l}M_{\beta}\left(\frac{1}{\sqrt{R_0}}+\frac{\sqrt{\lambda_k}}{\sqrt{a_0}}\right),\quad  M_{\beta}=\frac{1}{l^{\beta}\Gamma(2-\beta)}+1,\eqno(12)$$

\begin{proof}
  Since $u(x)\in H^{\beta}(0,l)$, by inequality (81) in \cite{Idczak}, we obtain the following estimate
$$\rvert I^{1-\beta}_{0,x}u(x) \rvert \leq \frac{1}{l}\rVert I^{1-\beta}_{0,x}u(x) \rVert_{L_1(0,l)}+\rVert D_{0,x}^{\beta}u\rVert_{L_1(0,l)} $$
$$\leq\frac{1}{l}\cdot\frac{l^{1-\beta}}{\Gamma(2-\beta)}\rVert u(x) \rVert_{L_1(0,l)}+\rVert D_{0,x}^{\beta}u\rVert_{L_1(0,l)}\leq \sqrt{l}M_{\beta}\left(\rVert u \rVert_{L_2(0,l)}+\rVert D_{0,x}^{\beta}u\rVert_{L_2(0,l)}\right). \eqno(13)$$
Considering $\rVert u\rVert_{L_2(0,l)}=1$ and the equality (9) in Section 3, the last equality yields the estimate (12).
\end{proof}

Now, we generalize the results in Proposition 1 and Proposition 2 to the case of a connected graph with a non-empty boundary. For this purpose, we use the well-known max-min principle.

According to the max-min principle,
$$\lambda_k=\sup_{A_{k-1}\subset\mathcal{H}^\beta_+(\mathcal{G})} \underset{\begin{smallmatrix}
 u\in A_{k-1}^{\bot } \\
 ||u||=1
\end{smallmatrix}}{\mathop{\inf }}\,
 \mathcal{R}(u)$$
where the supremum is taken over all $k-1$ dimensional linear subspaces $A_{k-1}$ of the space $\mathcal{H}_{+}^{\beta}(\mathcal{G})$.

Similarly as in (13), in the case of the graph, we obtain the following estimate
$$\rvert I^{1-\beta}_{0,x}u(x) \rvert \leq \frac{1}{l_i}\rVert I^{1-\beta}_{0,x}u(x) \rVert_{L_1(0,l_i)}+\rVert D_{0,x}^{\beta}u\rVert_{L_1(0,l_i)} $$
$$\leq\sum\limits_{i=1}^{n}\frac{1}{l_i}\cdot\frac{l_i^{1-\beta}}{\Gamma(2-\beta)}\rVert u(x) \rVert_{L_1(0,l_i)}+\sum\limits_{i=1}^{n}\rVert D_{0,x}^{\beta}u\rVert_{L_1(0,l_i)}\leq C_1 \rVert u(x) \rVert_{L_1(\mathcal{G})}+\rVert D_{0,x}^{\beta}u\rVert_{L_1(\mathcal{G})}$$
$$\leq C_2\left(\rVert u \rVert_{L_2(\mathcal{G})}+\rVert D_{0,x}^{\beta}u\rVert_{L_2(\mathcal{G})}\right). $$
where $C_1=\frac{1}{\Gamma(2-\beta) \underset{1\leq i \leq n}{\mathop{min}}\{l_i^\beta\}}+1$, $C_2=\sqrt {\rvert \mathcal{G} \rvert}C_1$. Taking into account that $\rVert u \rVert_{L_2(\mathcal{G})}=1$ and (9), we get
$$\rvert I_{0,x}^{1-\beta}u_k\rvert\leq C\left(1+\sqrt{\lambda_k}\right).$$

Further, we use notation $\lambda_k(\mathcal{G})$ for eigenvalues of the graph $\mathcal{G}$.

Together with $\mathcal{G}$, we consider the disconnected version of $\mathcal{G}$, which will be denoted by $[\mathcal{G}]$. In $[\mathcal{G}]$ we consider each vertex as a boundary vertex, i.e., in each edge $e_i=(0,l_i)$ we consider problem (10),(11).
The set of eigenvalues of the graph $[\mathcal{G}]$ is denoted by $\big\{\lambda_k[\mathcal{G}]\big\}_{k=1}^{\infty}$.
It is clear, that it can be obtained from $\bigcup\limits_{i=1}^n\big\{\lambda_k(0,l_i)\big\}$ by reordering it as non-decreasing sequence.

We put $\mathcal{H}_{+,Dir}^{\beta}[\mathcal{G}]=\bigoplus\limits_{i=1}^n\mathcal{H}_{+,Dir}^{\beta}(0,l_i)$.
From the inclusion $\mathcal{H}_{+,Dir}^{\beta}[\mathcal{G}]\subset\mathcal{H}_{+,Dir}^{\beta}(\mathcal{G})$ it follows, that
$$\left\{\underset{\underset{\rVert u \rVert=1}{\mathop{u\in A^{\bot}}}}{\mathop{inf}}\mathcal{R}(u)| A\subset \mathcal{H}^{\beta}_{+,Dir}[\mathcal{G}], \,dim(A)=k-1\right\}\subset \left\{\underset{\underset{\rVert u \rVert=1}{\mathop{u\in A}}}{\mathop{inf}}\mathcal{R}(u)| A\subset \mathcal{H}^{\beta}_{+,Dir}(\mathcal{G}),\, dim(A)=k-1\right\}.$$

So, using max-min principle, we get
$$\lambda_k[\mathcal{G}]\leq \lambda_k(\mathcal{G}).$$
From the last inequality, it follows, that
$$\sum\limits_{k=1}^{\infty}\frac{1}{\lambda_k(\mathcal{G})}\leq\sum\limits_{k=1}^{\infty}\frac{1}{\lambda_k[\mathcal{G}]}=\sum\limits_{i=1}^{n}\sum\limits_{k=1}^{\infty}\frac{1}{\lambda_k(0,l_i)}<+\infty.$$

Summarizing, we have

\textbf{ Theorem 2.} The series $\sum\limits_{k=1}^\infty\frac{1}{\lambda_k}$ is convergent and
$$\rvert u_k\rvert\leq C\left(1+\sqrt{\lambda_k}\right),$$
where the constant $C$ depends on $R_0$, $a_0$, $\beta$ and lengths of the edges of $\mathcal{G}$.

\section{Conclusion}
\,\,\,\, We considered fractional analog of the Sturm-Liouville problem on metric graph. Fractional Sturm-Lioville operator considered in our work consists combination of left Riemann-Liouville and right Caputo type fractional derivatives. Such combination of the operators are used to get simmetric and positive analogue of SLO. We proved that the considered operator has countable number eigenvalues, converging to infinity. Also, we investigated convergence of series $\sum\limits_{k=1}^\infty \frac{1}{\lambda_k}$ and get estimate for eigenfunctions. Our result improves similar result obtained in \cite{klimek2016},\cite{klimek2018}.\\

\textbf{Funding} The work was partially supported by the Ministry of Science and Higher Education of the Russian Federation. Grant ~075-02-2024-1447.

\end{document}